\documentclass[10pt]{amsart}

\usepackage[utf8]{inputenc}
\usepackage[english]{babel}
\usepackage{amsmath, amsfonts, amssymb, amscd, amsthm, mathabx}
\usepackage{color,soul}
\usepackage[dvipsnames]{xcolor}
\usepackage{enumerate}
\usepackage{hyperref}

\def\N{\mathcal{N}}
\def\indi{\boldsymbol{1}}
\def\Pr{\mathcal{P}}
\def\Qr{\mathcal{Q}}
\def\R{\mathbb{R}}
\def\Z{\mathbb{Z}}

\def\C{\mathbb{C}}
\def\S{\mathbb{S}}
\def\Sr{\mathcal{S}}
\def\K{\mathcal{K}}

\def\F{\mathcal{F}}

\def\T{\mathbb{T}}
 \newcommand\norm[1]{\|#1\|}

\newcommand{\funev}[2]{\left\langle{#1},\,{#2} \right\rangle}
\newcommand{\innerprodC}[2]{\left({#1}\mid {#2} \right)_{\C^N}}
\newcommand{\innerprodR}[2]{\left({#1}\mid{#2}\right)_{\R^N}}
\newcommand{\innerprodL}[2]{\left({#1}\mid{#2}\right)_{L^2}}

\newcommand{\norma}[1]{\left\|{#1}\right\|}

\DeclareMathOperator{\grado}{deg}

\DeclareMathOperator{\kernel}{ker}
\DeclareMathOperator{\img}{img}
\DeclareMathOperator{\coker}{coker}
\DeclareMathOperator{\dom}{dom}

\DeclareMathOperator{\real}{Re}
\DeclareMathOperator{\imag}{Im}
\renewcommand{\Re}{\real}
\renewcommand{\Im}{\imag}

\theoremstyle{plain}
\newtheorem{theorem}{Theorem}
\newtheorem*{theorem*}{Theorem}

\newtheorem*{lemma*}{Lemma}
\newtheorem{corollary}{Corollary}[theorem]

\theoremstyle{definition}

\theoremstyle{remark}

\title[SYSTEMS OF FUNCTIONAL-DIFFERENTIAL EQUATIONS]{PERIODIC SOLUTIONS FOR SYSTEMS OF FUNCTIONAL-DIFFERENTIAL
SEMILINEAR EQUATIONS AT RESONANCE}

\author{PABLO AMSTER, JULIAN EPSTEIN, ARTURO SANJUÁN}

\date{\today}

\begin{document}

\begin{abstract}
 Motivated by Lazer-Leach type results, we study the existence of periodic solutions for systems of functional-differential equations at resonance with an arbitrary even-dimensional kernel and  linear  deviating terms involving a general delay of the form $\int_0^{2\pi}u(t+s)\,d\lambda(s)$, where $\lambda$ is a finite regular signed measure. Our main technique shall be the Coincidence Degree Theorem due to Mawhin.
\end{abstract}


\maketitle

\section{Introduction}
\label{sec:introduccion}

The problem of finding periodic solutions to semilinear functional-differential equations at resonance has been widely studied. For instance, in the simplest situation of a scalar ordinary differential equation, nonlinear perturbations of harmonic oscillators of the form
\begin{equation}
  \label{eq:lazerleachscalar}
  x^{''}+m^2x + g(x)=p(t)
\end{equation}were considered by Lazer and Leach among others, see \cite{lazer1968} for $m=0$ and \cite{lazerleach} for strictly positive integer values of $m$. Equations like \eqref{eq:lazerleachscalar} are also known as Duffing equations at resonance.  In \cite{lazerleach}, periodic solutions to \eqref{eq:lazerleachscalar} were found assuming that $g$ is continuous, non-constant and bounded with limits at infinity, the forcing term $p$ is continuous, $2\pi$-periodic and   satisfies the well-known Lazer-Leach condition:
\begin{equation}
  \label{eq:lazerleachextended}
  \left|\int_0^{2\pi}p(t)e^{-imt}\,dt\right|<2\left|\lim_{s\to\infty}g(s)-\lim_{s\to-\infty}g(s)\right|.
\end{equation} 
The Lazer-Leach result was extended to systems. For instance in \cite{amsterdinapoli}, $2\pi$-periodic solutions $u:\R\to\R^N$,
for systems of size $N$ of the form
\begin{equation}
  \label{eq:lazerleachsys}
  u''(t)+m^2u(t)+g(u)=p(t),
\end{equation}
were found under appropriate conditions on the asymptotic behaviour on $g$ and a topological condition on the projection of $g(u)-p$ over the kernel of the linear operator $u\mapsto u'' +m^2 u$. The assumptions on the nonlinear term in the present work constitute an extension of the conditions in \cite{amsterdinapoli}, see  further discussion and Section~\ref{sec:general-conditions-on-N} for details.

Problems of the form \eqref{eq:lazerleachsys} may be generalized by introducing a deviating argument in order to obtain systems of functional-differential equations of the form
\begin{equation}
\label{eq:shiwangeqpre}
  u''(t)+m^2u(t)+g(u(t-\tau))=p(t)
\end{equation}with $\tau\in[0,2\pi)$. In \cite{shiwang}, periodic solutions to equations like \eqref{eq:shiwangeqpre} were found assuming a condition similar to
\eqref{eq:lazerleachextended}. In \cite{xianlongshunian}, these results were extended allowing unbounded sublinear nonlinearities.

A linear instance of \eqref{eq:shiwangeqpre} is the Gompertz equation with a time delay
\begin{equation}
  \label{eq:gompertzbasic}
  N'(t)=-\alpha N(t)\log(K/N(t-\tau))
\end{equation}where $N(t)$ is the density of a self-limiting cell population, 
$\alpha$ is a positive constant and $K>0$ is known as the environment carrying capacity. In fact, the change of variables $u=\log N$ yields  the linear equation
\begin{equation}
  \label{eq:gompertzlinear}
  u'(t)=\alpha u(t-\tau)+p
\end{equation}with $p=-\alpha \log K$. 

Extending this situation to the semilinear systems case, periodic solutions to systems of the form
\begin{equation}
  \label{eq:amsterkunapre}
  \begin{cases}
  u_1'(t)+a_1u_1(t)+b_1u_1(t-\tau_1)+g_1(u_1(t-\tau_1), u_2(t-\tau_2)) = p(t) \\
  u_2'(t)+a_2u_2(t)+b_2u_2(t-\tau_2)+g_2(t,{u_1}_t,{u_2}_t)=0
  \end{cases}
\end{equation}for periodic $p$ were found in \cite{amsterkuna}, assuming \eqref{eq:lazerleachextended} and the resonance condition $|a_1|<|b_1|$,  and $b_1e^{im\tau_1}=-a_1-im$ with $m:=\sqrt{b_1^2-a_1^2}\in\mathbb N$. Here the subscript notation $u_t$ stands for $u_t(\theta)=u(t+  \theta)$ and $g_2$ is a real valued continuous function defined on $\mathbb R\times C[-\tau_1,0]\times C[-\tau_2,0]$. 

It is important, in applications, to consider systems of semilinear functional-differential equations involving distributed delays such as $\int_{-\tau}^{0}u(t+s)\beta(s)\,ds$, where $\beta$ is some integrable function defined over $[-\tau,0)$.  Notice that, in the periodic case, $\int_{-\tau}^{0}u(t+s)\beta(s)\,ds=\int_0^{2\pi}u(t+s)\widetilde \beta(s)\,ds$, with $\widetilde \beta$ defined as the $2\pi$-periodic extension of 
the product $\beta\indi_{[-\tau,0]}:[-2\pi,0]\to \R$, where $\indi_{A}$ denotes the usual indicator function of a set $A$. Also notice that the discrete delay can be written as $u(t-\tau)=\int_{-\tau}^{0}u(t+s)\,d\delta_{-\tau}(s)=\int_0^{2\pi}u(t+s)\,d{\delta}_{2\pi-\tau}(s)$ where $\delta_{x}$ is the Dirac measure concentrated at $x$. Thus, discrete delays and distributed delays in periodic $N$-dimensional systems may be generalized altogether as expressions of the form $\int_0^{2\pi}d\lambda(s)u(s+t)$ with $\lambda$ an appropriate $N\times N$-matrix of regular signed measures. 

We are interested in finding $2\pi$-periodic solutions $u:\R\to\R^N$ to problems like \eqref{eq:lazerleachscalar}, \eqref{eq:lazerleachsys} and \eqref{eq:amsterkunapre} at resonance involving more general delays such as
\begin{equation}
\label{eq:maineqsimple}
  u'(t)+\int_0^{2\pi}d\lambda(s)u(s+t)+ g\left(\int_0^{2\pi}d\lambda(s)u(s+t)\right)=p(t).
\end{equation}Due to technical considerations we shall restrict this work to even dimensional kernels of the linear operator $L(u)(t):=u'(t)+\int_0^{2\pi}d\lambda(s)u(s+t)$. For situations including odd dimensional kernels, see \cite{lazer1968} for a very classical result or \cite{kaishiping}, for a more recent result.  A sufficient condition for $\kernel L$ to be even-dimensional and nontrivial consists in assuming that $\lambda$ is a skew-symmetric matrix of regular signed measures,
\begin{equation}
    \label{eq:nonsingularsimple}
    \det\int_0^{2\pi}d\lambda(s)\neq 0
\end{equation}and
\begin{equation}
\label{eq:nontrivialsimple}
    \det\left(ikI+\int_0^{2\pi}e^{iks}d\lambda(s)\right)=0
\end{equation}for some $k\in \mathbb Z$.   Conditions \eqref{eq:nonsingularsimple} and \eqref{eq:nontrivialsimple} can be derived straightforwardly by using Fourier transform. For every $k\in\mathbb Z$ satisfying \eqref{eq:nontrivialsimple},  
$\kernel L$ contains a nontrivial vector 
{of the form $\cos(kt)A + \sin(kt)B$ with $A,B\in\R^N$.}

We may naturally generalize the existence of limits at infinity in the scalar case to the multidimensional case by saying that  the radial limits $\lim g(sv)$ as $s\to\infty$ exist uniformly for  $v
\in \mathbb S^{N-1}$, the unit sphere of $\R^N$. However, the generalization of the Lazer-Leach  condition \eqref{eq:lazerleachextended} requires a technical discussion that can be found in Section \ref{sec:nirenbergconditions} below or in Section 3 of \cite{amsterdinapoli} and it is inspired by Nirenberg's condition in \cite{nirenberg}.  In this context, our general version of the Lazer-Leach conditions read as follows: 

{ 
\begin{equation}
\tag{LL1}\label{eq:lazerleachgen1}
    \mathcal P(g_w-p)\neq 0
\end{equation} 
for every $w\in \kernel L$ with $\|w\|_{L^2}=1$, where $g_w(t):=\lim_{s\to\infty}g(sw(t))$ a.e. and $\mathcal P$ is the orthogonal projection on $ \kernel L$, in the $L^2$-sense,
and 
\begin{equation}\tag{LL2}
\label{eq:lazerleachgen2}
\deg \gamma \neq 0    
\end{equation}
where $\gamma$ is the mapping given by
$\gamma(w):=\frac{\mathcal P(g_w-p)}{\|\mathcal P(g_w-p)\|_{L^2}}$ for   $w\in \kernel L$ with $\|w\|_{L^2}=1$, defined over the unit sphere of $\kernel L$ for the $L^2$ norm. } 

We are in condition of stating our main result, which shall be obtained as a direct consequence of a more general result  involving higher derivatives and different deviating terms for the linear and nonlinear terms (see Theorem \ref{mainresult} below). 

\begin{theorem}\label{th:1}
  Let $\lambda$ be a skew-symmetric square matrix of regular signed measures of  size $N$, let $g$ be a bounded-continuous vector field in $\R^N$ such that $\lim\limits_{s\to\infty}g(sv)$ exists uniformly for  $v\in \mathbb S^{N-1}$ and let $p$ be a continuous $2\pi$-periodic vector-valued function. Assume \eqref{eq:nonsingularsimple}, \eqref{eq:nontrivialsimple}, \eqref{eq:lazerleachgen1} and \eqref{eq:lazerleachgen2}, 
  then there exists a continuously differentiable $2\pi$-periodic solution of \eqref{eq:maineqsimple}.
\end{theorem}

In order to state and prove the general result, we introduce some notation and terminology in the following section.

\section{The general setting}
\label{sec:generalequation}

\subsection{Some notation}
We denote by $\T$ the quotient space $\R/2\pi\Z$. It is clear that the 
$2\pi$-periodic functions in $\R$ can be identified with those functions  defined on $\T$. We denote by $\mu$ the normalized Haar measure on $\T$, that is, $\mu$ is the unique Haar measure on $\T$ such that $\mu(\T)=1$. Thus, if $f:\T\to\C$ is a measurable function, then the integral of $f$ with respect to $\mu$ is given by $\int_\T
f\,d\mu=\frac{1}{2\pi}\int_0^{2\pi}f(s)\,ds$. For $k=1,2,\dots$ we denote by $C^k$ the space of all functions $u:\T\to\R^N$ such that $u$ has continuous derivatives up to order $k$. It is clear that $C^k$ is a Banach space with the norm $\norm{u}_{C^k}=|u|_\infty+|u^{(k)}|_\infty$, where $|u|_\infty=\max_t|u(t)|$. The space of continuous functions will be simply denoted by $C$ and the norm will be given by $|u|_\infty$. Let $m$ be a positive integer and, for $j=0,1,\dots,m$, let $A_j\in M_N(\R)$ be a $N\times N$ real matrix with $A_m$ non-singular. We denote by $P$ the matrix polynomial $P(x)=\sum_jA_jx^j$. If $\partial:C^1\subset C\to C$ is the differential operator defined  by $\partial u=u'$, then we denote by $P(\partial):C^m\subset C\to C$ the unbounded operator defined  by the formula $P(\partial)u=\sum_jA_ju^{(j)}$. 

Clearly, a system in terms of $P(\partial)u$ may be reduced to a higher rank system of first order equations. For this reason, it seems unnecessary
to consider systems of higher order derivatives. However, as we can check easily with some examples, to verify all the conditions stated below in the reduced first order system may be cumbersome.

As mentioned before, given $u\in C$, we shall denote by $u_t:C\to C$ the translation operator given by $u_t(s)=u(t+s)$.  Take for the moment $N=1$.  As we already noticed, for $\tau\in[0,2\pi)$, delays of the form $u(t-\tau)$ or $\int_{-\tau}^0u(t+s)\beta(s)\,ds$ may be regarded as $\funev{\Lambda}{u_t}$ where $\Lambda\in C^*$, the dual space of $C$. Thus, for any $N$ it is enough to consider $\Lambda\in L(C,\R^N)$ for a wide range of delays captured by the expression $\funev{\Lambda}{u_t}$. 

In this setting, we are interested in finding solutions $u\in C^1$ to the following ``abstract'' equation \begin{equation}
  \label{EQ}
  P(\partial)u+ \funev{\Lambda}{u_t}  + g(\funev{\Psi}{ u_t}) + h(t,u_t)=p(t),
\end{equation}
where $p\in C$, $\Lambda,\Psi\in L(C,\R^N)$ and $g\in C(\R^N,\R^N)$,  $h\in C(\T\times C, \R^N)$  are bounded functions. It is clear that equation \eqref{EQ}  is more general than \eqref{eq:maineqsimple}. We shall say that equation \eqref{EQ} is \textit{at resonance} if the linear term has nontrivial kernel. 

It is worth noticing that, for any $\Lambda\in L(C,\R^N)$ there exists an associated bounded linear operator  $\widetilde \Lambda: C\to C$ defined by the formula $\widetilde \Lambda u(t):=\funev{\Lambda  }{u_t}$.

\subsection{The linear and deviating term}

\label{sec:linear}

Let us define the unbounded linear operator $L:C^m\subset C\to C$ by the linear  part of equation \eqref{EQ}, that is: 
\begin{equation}
  \label{e:Ldef}
  L:= P(\partial) +  \widetilde \Lambda.
\end{equation}

Due to The Riesz Representation Theorem, to each $\Lambda\in L(C,\R^N)$ we may associate an $N\times N$ matrix $\lambda$ of finite regular signed measures which represents $\Lambda$ in the sense that $\funev{\Lambda}{u}=\int_\T d\lambda(t) 
u(t)$ for every $u\in C(\T,\R^N)$.

Let us denote by $L^2$ the complex Hilbert space of square-integrable functions defined on $\T$ with values in $\C^N$. The inner product of two functions $f,g\in L^2$ shall be denoted by $\innerprodL f g$ and is given by $\int_\T\innerprodC {f(t)}{g(t)}\,d\mu(t)$, where $\innerprodC x y=\sum_jx_j\overline{y_j}$. 
For $u\in L^2$ and $k\in \Z$, the Fourier Transform of $u$ at $k$ shall be given by $(\F u)(k)=\widehat u(k)=\int_\T u(t)e^{-ikt}\,dt\in \C^N$. It is clear that the set
$l^2:=\{\alpha:\Z\to \C^N: \sum_k|\alpha(k)|^2< \infty \}$, endowed with the inner product $(\alpha\mid \beta)_{l^2}=\sum_k\innerprodC{\alpha(k)}{\beta(k)}$ is also a complex Hilbert space and $\F:L^2\to l^2$ is a Hilbert space isomorphism. 
Also, we shall denote by $H^m$ the Sobolev space of $2\pi$-periodic functions in $L^2$ with weak derivatives in $L^2$ up to order $m$.

We define the Fourier transform of a matrix of complex measures $\lambda$ at $k$ as $\widehat \lambda (k) := \int_\T e^{-ikt}\,d\lambda(t)$, which is a complex square matrix of size $N$. 
The symmetric relations $\widehat{u}(-k)=\overline{\widehat{u}(k)}$ and $\widehat{\lambda}(-k)=\overline{\widehat{\lambda}(k)}$ are valid for every $k\in \Z$ because~\eqref{EQ} is an equation taking values in $\mathbb R^N$ for every $t$. 

As an example, consider $\funev{\Lambda}{u_t}=\widetilde \Lambda u(t)=Au(t)$ with $A=[a_{jl}]$ a real square $N$-matrix. 
The corresponding matrix $\lambda$ shall be represented by $[a_{jl}\boldsymbol\delta_0]$ with $\boldsymbol\delta_0$ the Dirac measure at 0. 
Moreover $\widehat{\lambda}(k)=A$ for every $k\in\Z$. Further examples will be discussed in Section~\ref{sec:examples}.

We shall use the Fourier transform in order to compute $\kernel L$ and $\img L$. Set $\widehat L:=\F \overline L\F^{-1}$ where $\overline L:H^m\subset L^2\to L^2$ is the natural extension of $L$. 
Then $\widehat L\widehat u(k)=L_k\widehat u (k)$ with $L_k:= P(ik) +\widehat{\lambda}(-k)\in M_N(\C)$. 

Let us illustrate this notation with an example. Take the linear part of the problem \eqref{eq:amsterkunapre}, then a direct computation shows that the matrices $L_k$ are given by
\[
L_k=\begin{bmatrix}
ik+a_1+b_1e^{-ik\tau_1} && 0 \\
0 && ik+a_2+b_2e^{-ik\tau_2} 
\end{bmatrix}.
\]
As mentioned, we are interested in equation \eqref{EQ} at resonance, namely, when $\kernel L \neq \{0\}$. Then, for some $k\in \Z$ the matrix $L_k$ should be singular and thus, the set $\K:=\{k\in \Z : \det L_k =0\}$ is non-empty.  
We claim that $\K$ is finite. Indeed, assume the contrary and for $k\in \K$, take $x_k\in \kernel L_k$ with $|x_k|=1$. On the one hand, $|\widehat{\lambda}(k)x_k|$ is uniformly bounded with respect to $k$ and, on the other hand $|P(ik_j)x_{k_j}|\to\infty$ for any ${k_j}\to\infty$. This is a contradiction and therefore $\K$ is finite.

Our assumptions on $L$ are the following: 
\begin{enumerate}[(L1)]
\item $\det L_0 \neq 0$.
\item $\kernel L_k=\kernel L_k^*$ for every $k\in \K$.
\end{enumerate}
Here, the ${}^*$ symbol stands for the Hilbert adjoint. Condition (L1), which is equivalent to say that $0\notin \K$, implies that $\K = \{\pm k_{1}\cdots \pm k_K\}$ with  $k_j\neq0$ for all $j$. 
When $\det L_0=0$, equation \eqref{EQ} leads to odd dimensional kernels, which require a different treatment and, as mentioned before, shall not be considered in this paper. 

Let us describe now the kernel of $L$. Denote with $\nu_k:=\dim\ker L_k$, and let $\{\Theta_{k,j}\in \C^N:j=1,\dots, \nu_k\}$ be an orthonormal basis for $\ker L_k$ and with $\Pr:C\to C$ the continuous projection
\begin{equation}
\label{eq:Pdef}
\Pr u (t) =  \sum_{k\in\K}\left\{\sum_{j=1}^{\nu_j} \left[
    \innerprodC{\widehat{u}(k)}{\Theta_{k,j}}\Theta_{k,j}\right] e^{ikt}\right\}.
\end{equation}
Then, a direct computation shows that $\kernel L = \img \Pr$ and $\dim\ker L=2\nu$, where $\nu:=\sum_{j=1}^K\nu_j$. For future reference we denote by
\begin{equation}
  \Sr :=\left\{u\in\kernel L : \norm{u}_{L^2}=1\right\}
\end{equation}
the unit sphere in $\kernel L$ in the sense of $L^2$. Since $\kernel L$ is finite-dimensional, using Fourier transform we notice that $\F (\mathcal S)=\sqrt 2\S^{2\nu-1}$.

The assumption (L2) 
compensates the lack of self-adjointness of the problem, i.e. 
the fact that not necessarily $\overline L= \overline L^\ast$. There are plenty of cases in which (L2) may occur. For instance, this
is the case if $\widehat{\lambda}(k)$ is skew-Hermitian and $P(x)$ has only odd powers or if $\widehat{\lambda}(k)$ is Hermitian and $P(x)$ has only even powers. Such is the case of the linear parts of  problems \eqref{eq:shiwangeqpre} and \eqref{eq:amsterkunapre}. It is clear that (L2) is equivalent to say that $\kernel \overline L=\kernel \overline L^*$.

Next, let us compute $\img L$. Take $\varphi\in C(\T,\R^N)$ such that $\widehat \varphi(k)\in\img L_k$ for $k\in \K$. Choose ${\phi}(k)\in \C^N$ such that $L_k\phi(k)=\widehat{\varphi}(k)$. We can define, by the Fourier representation
\begin{equation}\label{phirango} u(t)\sim \sum_{k\in \Z\setminus
\K}{(ikI+\widehat \lambda(-k))}^{-1}\widehat{\varphi}(k)e^{ikt} + \sum_{k\in \K}
\phi(k)e^{ikt}.
\end{equation}

It is not difficult to show, using again the fact that $\det A_m\neq 0$, that $u$ belongs to the space $ H^m $ which,  by Morrey's Inequality, is embedded in 
$C^{m-1,\frac12}$. Then, we can replace the $\sim$ sign in~\eqref{phirango} by $=$ everywhere in $t$. Using the structure
of~\eqref{EQ}, it is clear that $u\in C^1(\T,\R^N)$ and \begin{equation} \img L = \left\{\varphi\in C(\T,\R^N): \text{ for all
}k\in\K\text{, } \widehat\varphi(k)\in\img L_k={(\kernel L_k)}^\perp \right\}.
\end{equation}

Hence, $\img L$ is also a closed subspace in $C(\T,\R^N)$ and then a continuous projection $\mathcal Q:C\to \coker L$ may be given by $\mathcal Q=\mathcal P$ because $\coker L=\ker L$. Summing up, the unbounded operator
$L:C^m\subset C \to C$ with $\dom L =C^m$ is a zero-index Fredholm operator with the identity as the connecting isomorphism between the previous kernel and cokernel.

Now we need some compatibility assumptions between $\Lambda$ and the deviating argument $\Psi$. To this end, we shall also assume: 
\begin{enumerate}
\item[(L3)] $\det \widehat{\psi}(k)\neq 0$ for every $k\in \K$.
\item[(L4)] $\Theta_{k,j}$ is an eigenvector of $\widehat{\psi}(-k)$ for every $k\in\K$ and every
  $j=1,\dots,\nu_k$.
\end{enumerate}

Condition (L3) is a form of coerciveness in $\kernel L$. In fact, let us denote by $\widetilde \Psi:L^2\to L^2$ the bounded linear operator defined by $\widetilde \Psi (u )(t)= \funev{\Psi}{u_t}$. Then, (L3) is equivalent to say that there exists $c_\Psi>0$ such that, for every $z\in \kernel L$, $\norm{\widetilde \Psi z}_{L^2}^2=\sum_{k\in\K}|\widehat{\psi}(k)\widehat
z(k)|^2\geq c_\Psi^2 \norm{z}_{L^2}^2$. On the other hand, recalling \eqref{eq:Pdef}, condition (L4) ensures that $\widetilde\Psi$ is invariant on $\kernel L$. 

As an example to illustrate this situation, when $g$ does not depend on any delay, $\funev{\Psi}{u_t}=u(t)$ and (L3)--(L4) are trivially fulfilled and if $\Psi=\Lambda$, then condition (L4) is also fulfilled.

\subsection{The nonlinear term}
\label{sec:nonlinearterm}

Let us denote the nonlinear term or the Nemytskii operator in equation~\eqref{EQ} with $\N:C\to C$ and defined pointwise for every $t\in \T$ by
\begin{equation}
  \label{e:Ndef}
  \N u(t)= p(t)-g(\funev{\Psi}{u_t})-h(t,u_t),
\end{equation}
where $p\in C$. On the nonlinear perturbation $h(t,u_t)$, with $h: \T\times C\to \R$, we are going to assume only that $h$ is continuous and bounded, and we denote $|h|_\infty: = \sup\{|h(t,u_t)|:(t,u)\in\T\times C\}$.  On $g:\R^N\to \R^N$ we assume that is continuous and bounded, but we need some extra assumptions on its asymptotic behavior.

\subsection{A set of general conditions for the nonlinear term}
\label{sec:general-conditions-on-N}

Given $g$ as before, we define the following real valued function\begin{equation}
  \label{eq:guw-definition} g_{u,w}(t):=\limsup_{s\to\infty}\innerprodR{g(s\funev{\Psi}{u_t})}{w(t)}
\end{equation}
defined for $t\in\T$, $u\in C(\T)$ and $w\in \mathcal S$. Our assumptions on the nonlinear term $\N$, inspired by~\cite{amsterdinapoli}, are given by:

\begin{enumerate}[(N1)]
\item There is a family $\{(G_j,w_j)\}$ with $j=1,\dots,M$ such that $\{G_j\}$
  covers $\Sr $ and $w_j\in \Sr $. For every $j=1,\dots M$ and every $w\in U_j$, the map
  $g_{w,j}(t):=g_{w,w_j}(t)$ is upper semi-continuous at $w$ for almost every
  $t\in\T$. Namely, for almost every $t\in \T$, if $w_n\in C(\T)$ and
  $\norm{w_n-w}_\infty\to 0$, then
  \begin{equation*}
     \limsup_{w_n\to w}g_{w_n,j}(t)\leq
g_{w,j}(t).
  \end{equation*}
\item For each $w\in \Sr $, there exists $j\in\{1,\dots,M\}$ such that
  \begin{equation*}
  \overline{g_{w,j}} +|h|_\infty < {(p\mid w_j)}_{L^2}.
  \end{equation*}
  where $\overline{g_{w,j}}=\int_\T g_{w,j}\,d\mu$ is the average of the real
  valued function $g_{w,j}$ on the torus.
\item There exists $R_0>0$ such that, for every $R\geq R_0$,
  $\grado (\Gamma,B_R)\neq 0$, where $\grado(\Gamma,B_R)$ denotes the Brouwer degree at
  zero of $\Gamma$ in the open ball $\{\norma{w}_{L^2}<R\}$ and $\Gamma:\ker L\to \ker L$ is
  given by
  \begin{equation*} 
  {\Gamma}(w) := \Pr\N w.
\end{equation*}
\end{enumerate}

For $h\not\equiv 0$ the computation of $\Gamma$ may be difficult but, as we shall see later, for some systems with $h\not\equiv0$ that computation is unnecessary if $h$ does not interact with $\ker L$, in a way that shall be precised in Section~\ref{sec:amsterkunalike}.

Conditions (N1)--(N3) look cumbersome but they are quite general and, in some cases, they are easy to verify as shall be shown in the subsequent discussion.

\subsection{Nirenberg like conditions}
\label{sec:nirenbergconditions}

In this section, we shall present a set of conditions that imply (N1)--(N3) with $h\equiv 0$. Let us assume the following above mentioned radial condition:
\begin{enumerate}[(R1)]
\item The limit
  \begin{equation*}
    g_{v}:=\lim_{s\to\infty}g(sv)
  \end{equation*}
  exists uniformly with respect to $v\in \mathbb S^{N-1}$.
\end{enumerate}

Assuming (R1), for every $u\in C(\T)$ and for every $t\in\T$, we define
\begin{equation}
  \label{eq:radiallimitfunctions}
  g_u(t):=
  \begin{cases}
    \lim\limits_{s\to\infty}g(su(t)) &
    \text{if }u(t)\neq 0 \\
    g(0) & \text{in other case.}
  \end{cases}
\end{equation}

Fix $w\in\mathcal S$, we claim that for almost every $t\in\T$, (R1) implies that $g_w(t)$ is continuous with respect to $w$. Indeed, following \cite{amsterdinapoli}, take $w_n\in C(\T)$ such that $\norm{w_n-w}_\infty\to0$. It is easy to check that $T_w:=\{t\in\T:w(t)= 0\}$ is a null set. Fix $t\in\T\setminus T_w$. Because the relation $|w_n(t)-w(t)|<\frac12|w(t)|$ holds for all but finitely many values of $n$,
we may assume without loss of generality that $w_n(t)\neq0$ for every $n$. Let $\epsilon>0$, using (R1), we may choose $s_0>0$ such that $|g_{w_n}(t)-g(s_0w_n(t))|$ and $|g_{w}(t)-g(s_0w(t))|$ are smaller than $\epsilon/3$. Using the continuity of $g$, there exists $\delta_0>0$ such that for every $y\in\R^N$ with $|s_0w(t)-y|<\delta_0$, $|g(y)-g(s_0w(t))|<\epsilon/3$. There exists also $n_0\in\mathbb N$ such that for every $n\geq n_0$, $\|w_n-w\|_\infty<\delta_0/s_0$. Then we deduce that $|g_{w_n}(t)-g_{w}(t)|<\epsilon$, which proves our claim. Henceforth, (R1) implies (N1) for any family $\{(G_j,w_j)\}$ and  $g_{w,j}(t)=\innerprodR{g_w(t)}{w_j(t)}$.

Let us denote by $\widetilde\Gamma:\mathcal S\to \ker L$ the function given by 
\begin{equation}
  \label{eq:gamamdef}
{\widetilde{\Gamma}(w)}:= \Pr (g_w-p)
\end{equation}
and notice  that, for every $w\in\mathcal S$,
\begin{equation}
\label{eq:gamatilde}
  \lim_{s\to\infty}\norm{\widetilde\Gamma(w)- \Gamma(sw)}_{L^2}= 0.
\end{equation} 
We claim that the limit in \eqref{eq:gamatilde} is uniform with respect to $w\in\mathcal S$. In order to prove our claim, we shall need  the following tool from Approximation Theory, see \cite{nazarov}.
\begin{lemma*}[Nazarov-Turan Lemma]   Let $E\subseteq \T$ a $\mu$-measurable set and $q$ a polynomial with complex coefficients defined in $\T$. Then
  \begin{equation}
    \label{eq:nazarovturanineq}
    \norma{q}_{L^\infty}\leq \left(\frac{14}{\mu(E)}\right)^{m}\sup_{E}|q(t)|
  \end{equation}
where $m-1$ is the number of nonzero coefficients in $q$.
\end{lemma*}
Let $\epsilon>0$ and, for $w\in \mathcal S$, set $E(\epsilon):=\{t\in\T:|w(t)|<\epsilon\}$. Using Hölder inequality, norm equivalences in $\R^N$ and the Nazarov-Turan Lemma, we obtain, for for some $C_{NT}>0$:
\begin{equation}
  \label{eq:nazarov-turan-consequence}
  \mu(E_{\epsilon,w})\leq C_{NT}\epsilon^{\frac1{2\nu+1}}.
\end{equation}
 Note that the upper bound in \eqref{eq:nazarov-turan-consequence} is uniform with respect to $w\in\mathcal S$. Then, there exists $\widetilde\epsilon>0$ such that $\mu(E_{\widetilde\epsilon,w})<\epsilon^2/4$.

Using (R1), there exists $s_1>0$ depending only on $\epsilon$ such that  $|g_v-g(sv)|<\frac\epsilon2$ for every $v\in \S^{N-1}$. Take $s_0=\max\{s_1,s_1\widetilde\epsilon\}$. If $s\geq s_0$, then
\begin{equation}
  \begin{aligned}
    \label{eq:cuentalimiunifgamma}
    \|\widetilde\Gamma(w)-\Gamma(sw)\|_{L^2}^2 &= \int_ {\T\setminus
      E(\widetilde\epsilon)} |g(sw(t))-g_{\frac{w(t)}{|w(t)|}}|^2\,d\mu +
    \int_{E(\widetilde\epsilon)}|g(s(w(t))-g_{w(t)}|^2\,d\mu\\
    & <
    \mu(\T)\frac{\epsilon^2}4+\mu(E_{\widetilde\epsilon,w})|g|_\infty^2<\epsilon^2.
  \end{aligned}
\end{equation}
This completes the proof of the claim.

Next, we assume the following condition:
\begin{enumerate}
\item[(R2)] $\widetilde\Gamma(w)\neq 0$ for every $w\in \mathcal S$.
\end{enumerate}
If (R2) holds, then for every $w \in \mathcal S$, there exists a vector $w_j\in\mathcal S$ such that $(\widetilde\Gamma(w)\mid w_j)_{\R^N}<0$ in a neighborhood of $w$. Thus, compactness of $\mathcal S$ implies (N2).

Let us define $\gamma:\mathcal S\to \mathcal S$ by the formula
\begin{equation}
  \label{eq:gammadef}
  \gamma(w) = \frac{\widetilde\Gamma(w)}{\norm{\widetilde\Gamma(w)}_{L^2}}.
\end{equation}
Using Fourier transform, it is clear that $\gamma$ may be regarded as a mapping from $\mathbb S^{2\nu-1}$ into itself. Finally, we assume that 
\begin{enumerate}
\item[(R3)] $\deg \gamma \neq 0$.
\end{enumerate}
Let us prove that (R1)--(R3) imply (N3). Firstly, we claim that $(\Gamma,B_R)$ is admissible for $R\gg 1$ or, in other words, that the degree of $\Gamma$ over $B_R$ is defined. In fact, using (R2) and that the limit in \eqref{eq:gamatilde} is
uniform with respect to $w\in\mathcal S$, there is $R_0>0$ such that if $R\geq R_0$, then
\begin{equation}
\label{eq:gamma-admisible}
\|\Gamma(Rw)\|_{L^2}\geq \frac12\|\widetilde\Gamma(w)\|_{L^2}\geq
\frac12\min_{w\in\S}\|\widetilde\Gamma(w)\|_{L^2}>0.
\end{equation}
Then, inequality \eqref{eq:gamma-admisible} proves our claim. 
Further, observe that the  properties of the Brouwer degree imply that, for every $R\geq R_0$, $\grado (\Gamma,B_R)=\grado (\Gamma, B_{R_0})=\grado (\widetilde
\Gamma,\partial B_R)=\grado (\gamma)$.

Summing up, conditions (R1)--(R3) are a kind of Nirenberg conditions (see \cite{nirenberg}) and (R1)--(R3) imply (N1)--(N3) with $h\equiv 0$.

Let us discuss now the scalar case $N=1$. Assume (R1) which, as readily seen, just states the existence of the limits $g(\pm\infty)=\lim\limits_{s\to\pm\infty}g(s)$.
Let  $Lu=u''+m^2 u$ for some fixed $m\in \mathbb N$, then $\kernel L = \{A\cos(mt-\varphi):A\geq 0,\,\varphi\in\T\}$. That is, $\K=\{\pm m\}$. The classical Lazer-Leach
condition~\cite[p.~60]{lazerleach} reads
\begin{equation}
  \label{eq:classicll} |\widehat{p}(m)|<\frac 1\pi|(g(+\infty)-g(-\infty))|
\end{equation}
and implies (R2)--(R3). The proof of the last assertion may be found in~\cite{amsterdinapoli} but, for the reader's convenience we give an outline here. Let $w\in\Sr$, then
\begin{equation}g_w(t)=g(+\infty)\indi_{\{w(t)>0\}}(t)+g(-\infty)\indi_{\{w(t)<0\}}(t).\end{equation}
and henceforth, taking $w(t)=\sqrt 2\cos(mt-\varphi)\in\mathcal S$, we get: 
$\widehat{g_w}(m)=\frac1\pi(g(+\infty)-g(-\infty))e^{-i\varphi}$. This shows that if \eqref{eq:classicll} holds, then (R2) and (R3) are  verified.

\subsection{Limits at infinity}
\label{sec:limitsinfinity}

Let us denote by $\Sigma$ the set of all the $N$-tuples of signs $\sigma=(\pm,\dots,\pm)$ and, for $\sigma\in\Sigma$ and $x\in\R^N$, we shall write $\sigma x=(\sigma_1x_1,\dots, \sigma_Nx_N)$ where each $\sigma_j$ is $+$ or $-$. Similar to condition (R1), we
formulate the following condition for the $2^N$ different limits at infinity: 
\begin{enumerate}
\item [(R1')] For every $w\in\mathcal S$ and every $j\in\{1,\dots ,N\}$
  $\innerprodR{w(t)}{e_j}\not\equiv 0$ and, for every $\sigma\in\Sigma$,
    the  limits 
    \begin{equation}
      \label{eq:r1prime}
      \lim_{\substack{\sigma_jx_j\to\infty\\j=1,\dots N}}g(\sigma x) = g(\sigma)
    \end{equation}
exist.
\end{enumerate}
Assuming (R1'), each component of $w\in\mathcal S$, given by $\innerprodR{w(t)}{e_j}$ is a nontrivial null-average trigonometric polynomial. Thus, the set
\begin{equation}
  \label{eq:nullsetliminfinity}
  S_w:=\bigcup_{j=1}^NS_{w,j}:=\bigcup_{j=1}^N\left\{t\in\T: w_j(t)=0 \right\}
\end{equation}
is a null set (in fact, $S_w$ is a finite set). For $\sigma\in\Sigma$ let us define $w^{(\sigma)}:=\{t\in\T:\sigma_jw_j(t)>0,\,j=1,\dots,N\}$.

It is easy to prove that for $w\in\mathcal S$, the mapping $g_w(t)$ defined in \eqref{eq:radiallimitfunctions} is given by $g_w(t)=\sum_{\sigma}g(\sigma)\indi_{w^{(\sigma)}}$.  Moreover, a straightforward computation shows that, for every $t\in\T\setminus S_w$, $g_w(t)$ is
continuous with respect to $w$. In order to check that, in this context, (R1'), (R2) and (R3) imply (N1)--(N3), it is enough to follow the same argumentation line and techniques given in Section~\ref{sec:nirenbergconditions}.

\section{Main result and some consequences}
\label{sec:mainresults}

Under the notation given so far, our main result reads: 

\begin{theorem}
  \label{mainresult}
  Let $L$ and $\N$ be as in~\eqref{e:Ldef} and  \eqref{e:Ndef} respectively. Assume (L1)--(L4) and (N1)--(N3). Then the problem $Lu=\N u$ has a solution.
\end{theorem}

One immediate consequence of Theorem~\ref{mainresult} is  Theorem 3.1 of \cite{amsterdinapoli}, in which $C^2$ solutions to 
\begin{equation}
  \label{eq:lazerleachstructure}
  u''(t)+m^2u(t)+g(u)=p(t)
\end{equation}
with $m\in\mathbb N$ are found. Here $\funev{\Lambda}{u_t}=m^2u(t)$, $\funev{\Psi}{u_t}=u(t)$ and 
$\N u$, given by $g(u)-p$, satisfy (N1)--(N3).

As a consequence of Theorem \ref{mainresult}, taking into account the discussion in Section~\ref{sec:nirenbergconditions}, the following corollary is obtained.

\begin{corollary}
  \label{c:mainresult}
  Let $L$ be as in~\eqref{e:Ldef} and $\N$  be as in~\eqref{e:Ndef} with $h\equiv0$. Assume (L1)--(L3), (R2)--(R3) and either (R1) or (R1'). Then the problem $Lu=\N u$ has a solution.
\end{corollary}

It is clear that Corollary~\ref{c:mainresult} implies the classical Lazer-Leach result, see \cite{lazerleach}. 
\begin{theorem*}[Lazer-Leach]
  Let $g:\R\to\R$ be a bounded continuous non-constant function such that the
  limits $g(\pm\infty)$ exist. Let $p:\T\to\R$ be a continuous function such that
  \eqref{eq:classicll} condition holds. Then, the scalar equation \eqref{eq:lazerleachscalar}
  has a  $2\pi$-periodic $C^2$ solution.
\end{theorem*}

\section{Proof of the main result}
\label{sec:proofs}

For convenience, let us firstly state the standard topological tools that shall be  used in the  proof of our main result for the general case. 
Let $\mathbb X$ and $\mathbb Y$ be Banach spaces and $L: \dom L\subset \mathbb X\to \mathbb Y$ an unbounded zero-index Fredholm operator. Let $\Pr$ and $\Qr$ be projections on $\kernel L$ and $\coker L$ respectively and denote by $K:\img L\to \dom L\cap \ker \Pr$ the right inverse of $ L$ associated to $\Pr$. 
Let $\Phi:\coker L\to\ker L$ be a linear isomorphism. Let $ \N:\mathbb X\to \mathbb Y$ be a nonlinear operator such that $ K(I_{\mathbb Y}- \Qr)\N$ is completely continuous. With this notation in mind, our main tool is the following theorem due to Mawhin~\cite[Theorem~2.7]{mawhin}, which can be easily proven using the homotopy invariance property of the Leray-Shauder degree.

\begin{theorem}[Degree Continuation Theorem]\label{continuationtheorem}
  Let $\Omega$ be a bounded open set in $\mathbb X$ and suppose the following   conditions holds.
  \begin{enumerate}[({D}1)]
  \item $ Lx\neq \lambda  \N x$ for any
    $(x,\lambda)\in (\dom  L\cap \partial \Omega)\times (0,1)$.
  \item $ \Qr \N x\neq0$ for any
    $x\in \ker  L\cap \partial\Omega$.
  \item
    $\deg(\Phi \Qr\N,\Omega\cap \ker L)\neq 0$.
  \end{enumerate}
  Then the problem $ Lx= \N x$ has at least one solution in
  $\overline{\Omega}\cap \dom  L$.
\end{theorem}

In our case, for $L$ and $\N$ as in~\eqref{e:Ldef} and  in~\eqref{e:Ndef} respectively, the zero-index Fredholm property of $L$ is clear from the discussion in Section~\ref{sec:linear}. Moreover, using the norm $\|\,\cdot\,\|_{C^m}$ in $\ker \Pr$, it is clear that
the linear operator $\widetilde L=(I-\Qr)L(I-\Pr):\kernel \Pr \cap \dom L \to \kernel \Qr$ is continuous and bijective. In the present case $\mathcal P:C^1(\T)\to\ker L$ is defined, like before as in \eqref{eq:Pdef},
$\kernel \Pr \cap \dom L = (\ker L)\cap C^1 $, $\mathcal Q=\Pr$ and $\ker Q=\img L$. By the open mapping theorem, $\widetilde L^{-1}=K$ is also continuous and we have $|u'|_\infty\leq \kappa |Lu|_\infty$ for some constant $\kappa>0$ and every $u\in \ker P\cap \dom L$. The Arzelá-Ascoli Theorem implies that $K$ is compact. Then, the continuity and boundedness of $g$ and $h$ imply the compactness of $K(I-\mathcal Q)\N$.

Let us now check condition (D1) in Theorem~\ref{continuationtheorem}. By contradiction, suppose there exist $u_n\in C^1(\T,\R^N)$ and $\lambda_n\in (0,1)$ such that $\norm{u_n}_\infty\to\infty$ and $u_n$ satisfies $Lu_n=\lambda_n \N u_n$. If we write $u_n=v_n+z_n\in \ker P\oplus \img P$, then the Lyapunov-Schmidt reduction gives 

\begin{align} &v_n=\lambda_{n}K(I-\Qr)\N u_n\label{compacteq}\\ &\Qr\N
u_n=0.\label{resonanteq}
\end{align}

From~\eqref{compacteq}, $\norm{v_n}_\infty$ is bounded and, without loss of generality, we may assume that $\norm{z_n}_{L^2}\to\infty$. Define $\zeta_n(t):=\widetilde\Psi(z_{n})$, then it follows from the assumptions (L3)--(L4) on $\Psi$  that $\norm{\zeta_n}_{L^2}\geq c_\Psi\norma{z_n}_{L^2}\to\infty$ and $\zeta_n\in \ker L$. Setting  $\xi_n:=\zeta_n /\norm{\zeta_n}_{L^2}\in\Sr $ and going trough a subsequence if necessary, we may assume that $\xi_n\to \xi\in \Sr $ in
$L^2$. Let $w_j\in\Sr$ as in (N1), by (L2) we have that $\ker L = \ker L^*$ and \begin{equation}
  \label{e:pseudoselfadjoint} 0={(Lu_n\, |\, w_j)}_{L^2}={(u_n\mid L^*
w_j)}_{L^2} = {(\N u_n\mid w_j)}_{L^2}.
\end{equation}
Denoting $y_n=\xi_n+{\widetilde{\Psi}(v_{n})}/{\norm{\zeta_n}_{L^2}}$, equation~\eqref{e:pseudoselfadjoint} implies that
\begin{equation}\label{cuentaclave}
  \begin{aligned} {(p\mid w_j)}_{L^2}&=\limsup_{n\to \infty}
\left\{\int_\T \innerprodR{g(\funev{\Psi}{u_{nt}})}{w_j(t)}\,d\mu +
\int_\T \innerprodR{h(t,u_t)}{w_j(t)}\,d\mu\right\}\\ &\leq
\int_{\T}\limsup_{n\to \infty}\innerprodR{g(\funev{\Psi}{u_{nt}})}{
w_j(t)}\,d\mu + |h|_\infty\\ &\leq \int_{\T}\limsup_{n\to
\infty} \innerprodR{g(\norm{\zeta_n}_{L^2}y_n)}{w_j(t)}\,d\mu+ |h|_\infty \\
&\leq \int_{\T} g_{\xi,j}(t)\,d\mu + |h|_\infty=\overline{g_{\xi,j}} + |h|_\infty
  \end{aligned}
\end{equation}
which is contradiction with (N2). Thus condition (D1) holds.

In order to prove (D2), take $z_n\in\ker L$ with $\norm{z_n}_{L^2}\to\infty$ and assume that $\Qr \N z_n=\Pr\N z_n=0$. Then, reasoning as in \eqref{cuentaclave} we get again a contradiction with (L2). 

Finally, (D3) is a direct consequence of (N3) and the fact that, for every $w\in\ker L$, we have
\begin{equation}
  \Qr\N w  = \Pr\N w = \Gamma(w).
\end{equation}
This completes the proof of Theorem~\ref{mainresult}.

\section{Some examples and further discussion}
\label{sec:examples}

In this section we present some examples of equations and systems where Theorem~\ref{mainresult} and its Corollary~ \ref{c:mainresult} are applicable. For higher dimensional kernels with $\dim\ker L > 2$, explicit computations of $g_w$ and $\gamma$ may be difficult and extensive. For the sake of clarity, we shall   give here examples with $\dim\kernel L=2$.

\subsection{Duffing equation}
\label{sec:duffineq}
Fix $\tau\in[0,2\pi)$. In \cite{shiwang}, first-order systems of the form 
\begin{equation}
\label{eq:shiwangeq}
  u'(t)+Au(t)+g(u(t-\tau))+h(t,u_t)=p(t),
\end{equation}
where $A$ is an $N\times N$ matrix with pure imaginary eigenvalues are considered.  The Continuation Degree  Theorem~\ref{continuationtheorem} is employed as the  main technique and  sufficient conditions are obtained, expressed in a similar fashion of \eqref{eq:classicll}. Equation \eqref{eq:shiwangeq} is in the scope of this work under sufficient conditions given by (N1)--(N3) or some of its stronger versions. For instance, the authors consider the Duffing equation with a discrete delay as deviating argument, namely
\begin{equation}
  \label{eq:duffineq}
 u{''}(t)+m^2u(t)+g(u(t-\tau))=p(t).  
\end{equation}
Here, it is assumed that $p:\T\to\R$ is continuous and $g:\R\to\R$ is bounded and continuous such that the limits at infinity $g(\pm\infty)$ exist. 
In our setting, we are able to consider  equations with more general deviating terms; for instance, the equation
\begin{equation}
  \label{eq:duffingen}
   u{''}(t)+m^2u(t)+g(\funev{\Psi}{u_t})=p(t),
\end{equation}
where $\Psi\in C^*$ as long as $|\widehat\psi(m)|\neq0$, with $\psi$ the associated measure of $\Psi$. This is the case for a discrete delay since $\widehat\psi(m)=e^{im\tau}$, or a distributed delay $\funev{\Psi}{u_t}=\int_{-\tau}^0u(t+s)\beta(s)\,ds$ with $\beta\in L^1$ and $|\widehat\beta(m)|\neq0$. A direct computation shows that $\widehat {g_w}(m)=\frac1\pi e^{-i(\varphi+\theta_\Psi)}|g(+\infty)-g(-\infty)|$, where $w(t)=\sqrt2\cos(mt-\varphi)$ and $\theta_\Psi$ is fixed, given by $\cos\theta_\Psi=\frac{\Re\widehat \psi(m)}{|\widehat\psi(m))|}$ and $\sin\theta_\Psi=\frac{\Im\widehat\psi(m)}{|\widehat\psi(m)|}$. Then, the same classical Lazer-Leach condition \eqref{eq:classicll} is a sufficient condition which guarantees the existence of periodic solutions of \eqref{eq:duffingen}.

\subsection{A Gompertz system}
\label{sec:amsterkunalike}

It is worth noticing that it may happen, for some $h\not\equiv 0$, that $\innerprodR{h(t,u_t)}{\Theta_{k,j}}=0$ for all $t\in\T$, $u\in C $, $k\in\K$ and $j\in\{1,\dots,\nu_k\}$. In this case, the function $h$ does not get involved in the computation of the degree of $\Gamma$ or
$\gamma$. Let us exemplify this situation as follows.

Fix $\tau_1,\tau_2\in[0,2\pi)$. In \cite{amsterkuna}, periodic solutions of the system
\begin{equation}
  \label{eq:amsterkuna}
  \begin{cases}
    u_1'(t)+a_1u_1(t)+b_1u(t-\tau_1)+\widetilde{g}(u_1(t-\tau_1),u_2(t-\tau_2)) = p(t)\\
    u_2'(t)+a_2u_2(t)+b_2u(t-\tau_2)+h(t,u_t,v_t)=0
  \end{cases}
\end{equation}
were found assuming that 
\begin{enumerate}[(i)]
\item $\widetilde{g}$ is continuous and bounded, and the limits
  $$\widetilde{g}_{\inf}(\pm\infty):=\liminf\limits_{u\to\pm\infty} 
  \widetilde{g}(u,v), \qquad
  \widetilde{g}_{\sup}(\pm\infty):=\limsup\limits_{u\pm\infty} \widetilde{g}(u,v)$$
  exist uniformly with respect to $v$. \item $h:\T\times C\times C\to C$ is continuous and bounded. 
\item $a_1$, $b_1$ and $\tau_1$ satisfy the necessary and sufficient   condition for resonance at the first equation of  \eqref{eq:amsterkuna}. Namely $|a_1|<|b_1|$,   $m:=\sqrt{b_1^2-a_1^2}\in\mathbb N$ and $b_1e^{im\tau_1}=-a_1-im$.
\item $a_2$, $b_2$ and $\tau_2$ are such that $a_2+b_2\neq 0$ and do not satisfy the previous resonant condition for the second equation in
  \eqref{eq:amsterkuna}
\item
  $|\widehat p(m)| <
  \frac1\pi\max\{\widetilde{g}_{\inf}(+\infty)-\widetilde{g}_{\sup}(-\infty),
  \widetilde{g}_{\inf}(-\infty)-\widetilde{g}_{\sup}(+\infty)\}$.
\end{enumerate}

{Here, we shall assume that $\widetilde g_{\sup}=\widetilde g_{\inf}$ both at $+\infty$ and $-\infty$ and that \eqref{eq:classicll} holds for $\widetilde g$.} Thus, we may study systems like \eqref{eq:amsterkuna} with more general linear and deviating terms. Let us consider the system
\begin{equation}
  \label{eq:amsterkunagen}
  \begin{cases}
    u_1'(t)+\funev{\Lambda}{u_t} + \widetilde
    g(\funev{\Psi}{u_t},\funev{\Phi}{v_t})=p(t)\\
    u_2'(t)+\funev{\Phi}{v_t}+h(t,u_t,v_t)=0.
  \end{cases}
\end{equation}
with $\Lambda,\Psi,\Phi\in C^*$. Assume  there is a unique $m\in\mathbb N$ such that $im+\widehat\lambda(-m)=0$ and $im+\widehat\phi(k)\ne 0$ for all $k\in\Z$. Here $\widehat \lambda(k)$, $\widehat \phi(k)$ and $\widehat \psi(k)$ denote the Fourier transform at $k\in\Z$ of the signed measures associated $\Lambda$, $\Phi$ and $\Psi$ respectively. We need also to assume that $|\widehat\psi(m)|\neq 0$. Then, the linear term in \eqref{eq:amsterkunagen} satisfies (L1)--(L4) and $\mathcal S=\{(\sqrt2\cos(mt-\varphi),0)\}$. If we take
$g(u,v)=(\widetilde g(u),0)$, then it is clear that, for $w\in\mathcal S$, $g_w(t)$ is continuous for almost every $t\in\T$ and $\widehat{g_w}(m)=(\frac1\pi e^{i(\varphi+\theta_\Psi)}|\widetilde
g(+)-\widetilde g(-)|,0)$. Here $\theta_\Psi$ is fixed as in Section~\ref{sec:duffineq}. We conclude that system  \eqref{eq:amsterkunagen} has a solution.

\subsection{A weakly coupled system}
\label{sec:weakly-coupled}

Let $g(u)=(\widetilde g_1(u_1),\dots,\widetilde g_N(u_N))+ h(u)$ with each $g_i(\pm\infty)\in\R$ and $|h(u)|\to0$ uniformly as $|u|\to\infty$, then $g$ does not satisfy (R1) but it satisfies (N1).  In \cite{amsterdinapoli}, solutions of the system
\begin{equation}
  \label{eq:weaklycouple}
  u_j''(t)+m^2u_j+\widetilde g_j(u_j)+h_j(u)=p_j(t)\quad j=1,\dots,N.
\end{equation}
were found under the assumption
\begin{equation}
  \label{eq:clasicllcom}
  |\widehat p(m)|<\frac1\pi|\widetilde g_j(+\infty)-\widetilde g_j(-\infty)|\quad j=1,\dots, N.
\end{equation}

Let us consider the more general system
\begin{equation}
  \label{eq:weaklycouplegen}
  u_j''(t)+ \funev{\Lambda}{u_{j,t}} +\widetilde
    g_j(\funev{\Psi}{u_{jt}})+h_j(u)=p_j(t)\quad j=1,\dots,N.
\end{equation}
where $\Lambda$ and $\Psi$ are as in Section~\ref{sec:amsterkunalike}. If we assume \eqref{eq:clasicllcom}, then system \eqref{eq:weaklycouplegen} has a solution. This can be shown by noticing that for $R\gg 1$ there exists an admissible homotopy between $\Gamma$ and $T(w)=\widehat p(m)-\frac1\pi e^{-i(\varphi+\theta_\Psi)}[|\widetilde
g_j(+\infty)-\widetilde g_j(-\infty)|]_{j=1}^N$ representing $w(t)=\sqrt2\cos(mt-\varphi-\theta_\Psi)$.

\subsection{Distributed Delays}

As an example of an equation with non-trivial kernel and distributed delay, let us consider the problem 
\begin{equation}
\label{eq:distributed}
u'(t)+\alpha \int_{-\tau}^0 \beta(\theta)u(\theta+t)\,d\theta + g\left(\int_{-\tau}^0 \beta(\theta)u(\theta+t)\,d\theta \right) = p(t).
\end{equation}
with $\tau=\pi/m$ for some positive integer $m$. If we take $\alpha = m/2$ and $\beta(s)$ as the uniform distribution $\beta(s)=\boldsymbol{1}_{[-\pi/m,0]}$ or if we take $\alpha =4m/\pi$ and $\beta(s)$ as the density function given by $\frac{-m}2\sin(mx)\boldsymbol{1}_{[-\pi/m,0]}$, then we obtain non-trivial kernels, conditions (R1) to (R3) are also satisfied and hence there exist periodic solutions to \eqref{eq:distributed}

\subsection{Higher dimensional kernels}

Despite the verification of (R2) and (R3) may be cumbersome and might require a numerical approach, we may also consider equations with higher dimensional kernels such as the following beam-like equation:
\begin{equation}
    \label{eq:beamlike}
    u^{(iv)}(t)+(m_1^2+m_2^2)u''(t)+m_1^2m_2^2u(t) +g(\funev{\Psi}{u_t})=p(t)
\end{equation}
where $m_1,m_2\in\mathbb Z_{>0}$. 

\bibliographystyle{plain}


\end{document}